\documentclass[12pt,a4paper]{article}
\usepackage{amsmath,amstext,amssymb,amscd}

\begin{document}

\title{The purity of set-systems related to Grassmann necklaces}
      \author{V.I.~Danilov\thanks{Central Institute of Economics and
Mathematics of the RAS, 47, Nakhimovskii Prospect, 117418 Moscow, Russia;
email: danilov@cemi.rssi.ru},
 A.V.~Karzanov\thanks{Institute for System
Analysis of the RAS, 9, Prospect 60 Let Oktyabrya, 117312 Moscow, Russia;
email: sasha@cs.isa.ru.},
 and G.A.~Koshevoy\thanks{Central Institute of Economics
and Mathematics of the RAS, 47, Nakhimovskii Prospect, 117418 Moscow,
Russia; email: koshevoy@cemi.rssi.ru}}

\date{}
      \maketitle

\section{Introduction}
Studying the problem of quasicommuting quantum minors, Leclerc and Zelevinsky
\cite{LZ} introduced the notion of weakly separated sets in $[n]:=\{1,\ldots,
n\}$. Moreover, they raised several conjectures on the purity for this
symmetric relation, in particular, on the Boolean cube $2^{[n]}$ (or the
max-clique purity of the graph on $2^{[n]}$ generated by this relation). Recall
that a finite graph $G$ is {\em pure} if all (inclusion-wise) maximal cliques
in it are of the same cardinality. In \cite{max} we proved these purity
conjectures for the Boolean cube $2^{[n]}$, the discrete Grassmanian
${[n]\choose r}$, and some other set-systems. In \cite{OPS} the purity was
proved for weakly separated collections inside a positroid which contain a
Grassmann necklace $\mathcal N$ defining the positroid. We denote such
set-systems as $\mathcal{I}nt(\mathcal N )$; they are special collections of
sets in the discrete Grassmannian. The discrete Grassmannian itself is such a
collection for the largest necklace.

In this paper we give an alternative (and shorter) proof of the purity of
$\mathcal{I}nt(\mathcal N )$ and present a stronger result. More precisely, we
introduce a set-system  $\mathcal{O}ut(\mathcal N )$ complementary to
$\mathcal{I}nt(\mathcal N )$, in a sense, and establish its purity. Moreover,
we prove (Theorem~3) that these two set-systems are weakly separated from each
other. In the proof of this theorem, we use a technique of plabic tilings from
\cite{OPS}. As a consequence of Theorem~3, we obtain the purity of set-systems
related to pairs of weakly separated necklaces (Proposition 4 and Corollaries 1
and 2). Finally, we raise a conjecture on the purity of both the interior and
exterior of a generalized necklace. Our study of some other pure set-systems is
given in~\cite{pur2}.

           \section{Preliminaries}

For a natural number $n$, we denote by ${[n]\choose r}$ the set of
$r$-element subsets in  $[n]:=\{1,\ldots,n\}$ (the discrete Grassmanian).
Subsets of ${[n]\choose r}$ are called {\em (set-)systems} and we use
calligraphic letters for them.

It will be convenient for us to think of $[n]$ as being $\mathbb{Z}$ modulo
$n$. We consider the cyclically shifted orders $<_i$ on $[n]$, $i=1,\ldots,n$,
defined by $              i <_i (i+1) <_i\ldots <_i n <_i 1 <_i \ldots <_i
(i-1)$. A sequence $i_1, \ldots,i_k$ is called {\em cyclically ordered} if
$i_1<_i i_2<_i \ldots <_i i_k$ for some $i$.

We denote by  $\ll _i$ the following binary relation on ${[n]\choose r}$. For
two sets $X$ and $Y$ of cardinality $r$, we write $X \ll _i Y$ if  for any
$x\in X-Y$ and $y\in Y-X$, one holds $x <_i y$. \medskip

{\bf Definition.} Two subsets $X,Y\subset[n]$ of the same
cardinality\footnote{The definition of weak separability can be given for
arbitrary subsets in $[n]$; see \cite{LZ,Sc,max,OPS}. But in this paper we
deal only with the above-mentioned case.} are called weakly {\em
separated} (denoted as $X \| Y$) if $X\ll _j Y$ holds for some $j\in
[n]$.\medskip

In general, the relation $\ll_i$ is
not transitive. Nevertheless, the following assertion is valid.
\medskip

{\bf Lemma 1.} \cite[Lemma 3.6]{LZ} {\em Let $X\ll_i Y \ll_i Z$, where
$X,Y,Z$ have the same cardinality, and $X$ and $Z$ are weakly separated. Then $X\ll_i Z$.}\medskip

The notion of weak separation has proved its usefulness in the study of
Pl\"ucker coordinates on Grassmannians. Since we never deal with the
strong separation in this paper, we will use the term `separation' instead
of `weak separation' for short.

It is easy to see that  $X\ll _i Y$ for some $i\in [n]$ if and only if $Y\ll _j
X$ for some $j$. Therefore, the separation relation $\|$ on ${[n]\choose r}$ is
symmetric and reflexive. We say that two set-systems $\mathcal X $ and
$\mathcal Y $ from ${[n]\choose r}$ are \emph{separated from each other} (and
write $\mathcal X \|\mathcal Y $) if $X\| Y$ for any $X\in \mathcal X$ and
$Y\in \mathcal Y $. A system $\mathcal X $ is called {\em separated} if
$\mathcal X \|\mathcal X$. A system $\mathcal D \subset {[n]\choose r}$ is
called {\em pure} if all maximal separated subsystems in  $\mathcal D $ are of
the same size; this size is called the \emph{rank} of $\mathcal D $ and denoted
by $rk(\mathcal D )$.

We will essentially use the following important fact. \medskip

{\bf Theorem 1.} {\em The Grassmannian ${[n]\choose r}$ is a pure system of
rank  $r(n-r)+1$.}\medskip

This assertion was conjectured in \cite{LZ, Sc} and answered affirmatively
in \cite{max}. In fact,~\cite{max} proved the purity of the Boolean cube
$2^{[n]}$, and the above theorem follows from the argument of Leclerc and
Zelevinsky in~\cite{LZ} that the purity of the Boolean $n$-cube would
imply the purity of the Grassmannians ${[n]\choose r}$.

In  \cite[theorem 4.7]{OPS} the purity was shown for some systems of more
general character in ${[n]\choose r}$; they are produced from the so-called
Grassmann necklaces.

In the next section we recall necessary definitions. Throughout the paper,
symbol $\subset$ stands for non-strict inclusion (admitting equality).

           \section{Necklaces and related set-systems}

{\bf Definition.} \cite{Po} A (Grassmann) {\em necklace} in ${[n]\choose
r}$ is a family $\mathcal N =(N_1,\ldots,N_n)$ of sets from ${[n]\choose
r}$ such that $N_{i+1}$ contains $N_i-\{i\}$ for each $i$ (hereinafter the
indices are taken modulo $n$).\medskip

In particular, if $i\notin  N_i$ then $N_{i+1}=N_i$ and $i\notin N_j$ for all
$j$. We will assume for simplicity (see Remark 2 below) that this is not the
case, and that any $i\in [n]$ satisfies $i\in N_i$.

The necklaces are closely related to permutations on $[n]$. The set
$N_{i+1}$ is obtained from $N_i$ by deleting $i$ and adding some element
$\pi (i)$ (which may coincide with $i$). Thus, the necklace  $\mathcal N $
defines the corresponding map $\pi : [n] \to [n]$. This $\pi $ is
bijective. (Indeed, suppose that some element $j$ is not used. Then it
occurs either in none $N_i$ (which contradicts $j\in N_j$) or in all $N_i$
(yielding $j=\pi (j)$).) Therefore, $\pi $ is indeed a permutation on
$[n]$.

Conversely, let $\pi $ be a permutation on $[n]$. We can associate to it
the following family of sets $\mathcal N =\mathcal N _\pi
=(N_1,\ldots,N_n)$ by the rule
      $$
                         N_i=\{j\in [n], j\le _i \pi ^{-1}(j)\}.
      $$
It is easy to see that $\mathcal N$ is a necklace in ${[n] \choose r}$,
where the number  $r$ is defined to be the `average clockwise rotation' by
$\pi$ of the elements of $[n]$.\medskip

{\bf Example 1.} Let a permutation $\pi $ send every $i$ to $i+r$
(`rotation' by $r$ positions). Then $N_i=\{i,i+1,\ldots,i+r-1\}=[i,i+r)$
is a cyclic interval of length $r$ beginning at $i$. The corresponding
necklace is called the {\em largest} one; this terminology will be
justified later.\medskip

An important property of necklaces is given in the following
\medskip

{\bf  Lemma 2}. (\cite[Lemma 4.4]{OPS}) For all $i$ and $j$, one holds
$N_i-N_j\subset [i,j)=\{i,i+1,\ldots,j-1\}$.\medskip

Symmetrically, $N_j-N_i\subset [j,i)$. As a corollary, we obtain that $N_i\ll
_i N_j$ for any $i$ and $j$. In particular, all sets in a necklace $\mathcal N
$ are separated from each other.

For a necklace $\mathcal N$, let us call  the {\em interior} of $\mathcal N $
the following set-system
\[
\mathcal{I}nt(\mathcal N )=\{X\in {[n]\choose r}, \ N_i\ll _i X
      \text{ for every  }i\}.
\]
Obviously, $\mathcal N \subset \mathcal{I}nt(\mathcal N )$ and $\mathcal N
\|\mathcal{I}nt(\mathcal N )$.\medskip

{\bf A supplement to Example 1.} Let $\mathcal N $ be the largest necklace
consisting of cyclic intervals (see  Example 1). Since $[i,i+r)\ll _i      X$
for any $r$-element set $X$, we obtain that the interior of $\mathcal N $ is
the discrete Grassmanian, $\mathcal{I}nt(\mathcal N )={[n]\choose r}$. This
justifies the term `largest': this necklace has the largest interior.

Theorem 1 asserts that the interior of the largest necklace is a pure system.
This is generalized as follows.
\medskip

{\bf Theorem 2.} {\em For every Grassmann necklace $\mathcal N $, the set-system
$\mathcal{I}nt(\mathcal N )$ is pure.}\medskip

{\bf Remark 1}. This result is obtained in \cite{OPS}. Strictly speaking,
\cite{OPS} considered another system $Pos(\mathcal N)$, a
\emph{positroid}, and the purity is proved only for weakly separated
systems $\mathcal C \subset Pos(\mathcal N )$ which contain $\mathcal N $,
$\mathcal N\subset \mathcal C$. It is rather easy to show that such
systems are exactly weakly separated systems in $\mathcal{I}nt(\mathcal N
)$. Therefore, Theorem 2 is equivalent to Theorem 4.7 in \cite{OPS}. A
question on the purity of the positroid $Pos(\mathcal      N )$ (without
the additional condition $\mathcal N \subset \mathcal C $) is
open.\medskip

{\bf Remark 2}. Suppose that $i\notin N_i$ for some $i$. Then $i\notin X$ for
every $X\in \mathcal{I}nt(\mathcal N )$. Indeed, supposing $i\in X$, we obtain
a contradiction to $N_i \ll_i X$. Deleting such dummy $i$'s, we may assume that
$i\in N_i$ for any $i\in [n]$.

We give an alternative proof of Theorem 2 in the next section.

          \section{Alignments and extensions of necklaces}

To prove Theorem 2, it is convenient to consider another description for
the system $\mathcal{I}nt(\mathcal N )$, given in terms of alignments of
the permutation $\pi =\pi (\mathcal N )$. We use the notion of an
alignment introduced by Postnikov  \cite{Po}. Let $\pi $ be a permutation
of $[n]$. A pair $(i,j)$ is said to be an \emph{alignment} for $\pi $ (and
denoted by $i \Rightarrow_\pi j$) if the quadruple  $\pi ^{-1}(i), i,j,
\pi ^{-1}(j)$ occurs in this cyclical order (the case $j=\pi (j)$ is
admitted, whereas $i=\pi (i)$ is not). Roughly speaking, the `arrows'
entering $i$ and $j$, go parallel (do not cross) and in the same
direction. See the picture.

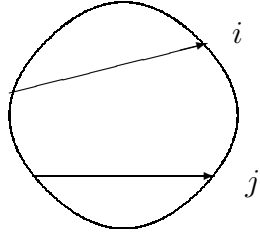
\unitlength=1mm \special{em:linewidth 0.4pt} \linethickness{0.4pt}
\begin{picture}(80.00,40.00)(-10,10)
\bezier{112}(50.00,20.00)(60.00,10.00)(70.00,20.00)
\bezier{112}(50.00,40.00)(40.00,30.00)(50.00,20.00)
\bezier{112}(50.00,40.00)(60.00,50.00)(70.00,40.00)
\bezier{112}(70.00,40.00)(80.00,30.00)(70.00,20.00)
\put(45.00,33.00){\vector(4,1){26.00}}
\put(48.00,22.00){\vector(1,0){24.00}}
\put(75.00,41.00){\makebox(0,0)[cc]{$i$}}
\put(77.00,21.00){\makebox(0,0)[cc]{$j$}}
\end{picture}

Notation $i \Rightarrow_\pi j$ for the alignment is justified as follows. Let
$Y\in\mathcal{I}nt(\mathcal N )$. If $i\in Y$ satisfies the relation
$i\Rightarrow_\pi j$, then $j\in Y$. We call this property of $Y$ the {\em
$\pi$-chamberness}. Indeed, without loss of generality, one may assume that
$i=\pi (1)$. Then $i$ does not belong to $\mathcal N _1$, whereas $j\in
\mathcal N _1$. Now suppose that  $i\in Y$ and $j\notin Y$. Then $i\in
Y-\mathcal N_1$ and $j\in \mathcal N _1-Y$. Due to the relation $N_1\ll _1 Y$
we obtain $j<_1 i$, which contradicts $i<j$.

The converse property takes place as well.\medskip

{\bf Proposition 1.} {\em For a set $Y\subset [n] $ of size $r$, the
following statements are equivalent:

           1) $Y\in \mathcal{I}nt(\mathcal N (\pi ))$,

           2) $Y$ is $\pi$-chamber set.}\medskip

The implication  $1) \Rightarrow  2)$ has been proved. To see the
implication $2) \Rightarrow  1)$, we show that 2) implies $N_i\ll _i Y$
for any $i$. Without loss of generality we may assume that $i=1$; so we
have to prove that $N_1 \ll_1 Y$. Suppose this is not so, i.e., there
exist $j\in \mathcal N _1-Y$ and       $i\in Y-\mathcal N _1$ such that
$i<j$. Then $j\in \mathcal N _1$ means that $\pi ^{-1}(j)>j$; and $i\notin
\mathcal N _1$ means that $\pi ^{-1}(i)<i$. This together with the
inequality $i<j$ means that the pair $(i,j)$ is an alignment. But then the
chamberness of $Y$ implies that $j\in Y$ (since  $i\in Y$). A
contradiction. \hfill$\Box$\medskip

In what follows we write $\mathcal{I}nt(\pi )$ for $\mathcal{I}nt(\mathcal N
)$.

We prove Theorem 2 by induction on the number of alignments of the
permutation $\pi $ corresponding to a necklace $\mathcal N$.
\smallskip

1. {\em A base of induction}: there are no alignments. In this case the
permutation $\pi $ sends each $i$ to  $i+r$. Indeed, let $\pi$ send  $i$ to
$i+k(i)$, $0< k(i)\le n$. Choose $i$ with $k(i)$ minimum. Then in case $\pi
(i-1)>\pi (i)$, we have $i-1\Rightarrow_\pi i$. This is impossible; so $\pi
(i-1)<\pi (i)$. Hence, $k(i-1)\le k(i)$. The minimality of $k(i)$ gives
$k(i-1)=k(i)$. Repeating this procedure, we obtain that $k(\cdot)$ is constant
(and equal to $r$).

Hence, the necklace with no alignments is the largest necklace and the
proposition follows from Theorem 1.\smallskip

2. {\em A step of induction.} Suppose that the permutation $\pi$ has an
alignment. Then there exists a `simple' alignment $i \Rightarrow_\pi j$,
in the sense that $\pi ^{-1}(i)$ and $\pi ^{-1}(j)$ are (cyclically)
consecutive numbers. Without loss of generality, we may assume that the
first number is 1 and the second one is $n$, so that $i=\pi (1)$ and
$j=\pi (n)$.

Now we consider the permutation $\pi '$ which coincides with $\pi $
everywhere except for the elements  $1$ and $n$. More precisely,
$\pi'(1)=j$ and $\pi '(n)=i$. If for the permutation $\pi$, the arrows
going from $1$ and from $n$ do not cross (and therefore give a simple
alignment $i\Rightarrow j$), then similar arrows for $\pi'$ do cross (and
the alignment $i\Rightarrow j$ vanishes). All other alignments preserve.
Thus, the set of alignments for $\pi '$ is obtained from that of $\pi$ by
deleting one alignment $i\Rightarrow_\pi j$. By induction the set-system
$\mathcal{I}nt(\pi ')$ is pure (and contains $\mathcal{I}nt(\pi )$, as
follows from  Proposition 1). Now Theorem 2 follows from the following
\medskip

{\bf Proposition 2.} {\em  Let $\pi$ and $\pi'$ be as above, let $X$ be a
set in $\mathcal{I}nt(\pi ')$ which is separated from $N_1$ and such that
$X\neq N'_1$. Then $X\in      \mathcal{I}nt(\pi )$.}\medskip

Indeed, let  $\mathcal{C}$ be a maximal separated subsystem in
$\mathcal{I}nt(\pi )$. Then the system $\mathcal{C}\cup\{N_1'\}$ is
contained in  $\mathcal{I}nt(\pi ')$ and is weakly separated. We assert
that it is a maximal separated system in $\mathcal{I}nt(\pi ')$. For if
this is not so, we can add some $X$ to this system. Then, due to
Proposition 2, $X$ belongs to $\mathcal{I}nt(\pi)$, which contradicts the
maximality of $\mathcal{C}$ in $\mathcal{I}nt(\pi)$. Thus,
$\mathcal{I}nt(\pi)$ is pure and the rank of $\mathcal{I}nt(\pi )$ is less
by $1$ than the rank of $\mathcal{I}nt(\pi ')$. By the  induction, we
conclude that the rank of $\mathcal{I}nt(\pi )$ is equals to $k(n-k)+1$
minus the number of alignments for $\pi $. This gives Theorem 2. \medskip

{\em Proof of Proposition 2}. Let $X$ be as in Proposition 2. We assert that
$X$ belongs to $\mathcal{I}nt(\pi )$. Suppose, for a contradiction, that $X$ is
not a $\pi$-chamber set. Since $X$ is a $\pi'$-chamber set and $\pi $ has
exactly one additional alignment $i \Rightarrow j$ compared with $\pi '$, we
have $i\in X$ and $j\notin  X$. The set $N'_1$ also contains $i$ but not $j$.
(Recall that $N'_1$ differs from $N_1$ by swapping the roles of $i$ and $j$:
$N_1$ contains $j$ and does not contain $i$.) Our aim is to prove that $X$
coincides with $N'_1$.

We have $N_n\ll _n X$. This means that any element of $X-N_n$ is greater
by $>_n$ than any element of $N_n-X$. Since $i$ belongs to  $X$ and does
not belong to $N_n$ (as $i$ appears only in $N_2$), we conclude that any
element of  $N_n-X$ is $<_n i$. Hence, besides $n$, any element of $N_n-X$
is $<i$. In other words, within the interval $I=(i,n)$ we have the
inclusion $N_n\subset X$. In this $I$ the sets $N_n$ and $N'_1$ coincide;
so within $I$ we have the inclusion $N'_1\subset X$. Since $n\notin N'_1$
(as $n$ is replaced by $i$ under changing $N'_n=N_n$ to $N'_1$), the set
$N'_1$ is contained in $X$ within $(i,n]$.

Similarly, using the relation $N_2\ll _2 X$, we obtain that $X\subset
N'_1$ on the interval $[1,j)$. In particular, within $(i,j)$ (and even
within $[i,j]$) the sets  $X$ and $N'_1$ coincide.

If the inclusion  $X\cap [1,i)\subset N'_1\cap [1,i)$ is strict, then the
inclusion  $N'_1\cap (j,n]\subset X\cap (j,n]$ is also strict. Hence, there are
an element $i'<i$ belonging to $N'_1-X$ and an element $j'>j$ belonging to
$X-N'_1$. Since $N'_1$ and $N_1$ coincide outside $\{i,j\}$, the element $i'$
belongs to $N_1-X$, and $j'$ belongs to $X-N_1$. Recall also that $i\in X-N_1$
and  $j\in N_1-X$. These relations together with the inequalities $i'<i<j<j'$
imply that the sets $N_1$ and $X $ are not weakly separated. This contradiction
completes the proof of Proposition 2. \hfill$\Box$

          \section{Exterior of a necklace}

In this section we show the purity of the so-called exterior  of a necklace.
Denote by $\mathcal S(\mathcal N)$ the system of sets weakly separated from the
necklace $\mathcal N$:
\[
\mathcal S(\mathcal N):=\{X\in {[n]\choose r}, ~X\|N_i,\,\forall\,i\in
[n]\}.
\]
We know that $\mathcal{I}nt(\mathcal N )$ is a subset of $\mathcal S(\mathcal
N)$. The {\em exterior } of a necklace  $\mathcal N $, denoted as
$\mathcal{O}ut(\mathcal N )$, is the complement to $\mathcal{I}nt(\mathcal N )$
in  $\mathcal S(\mathcal N)$, that is $\mathcal{O}ut(\mathcal N )= \mathcal
S(\mathcal N)\setminus \mathcal{I}nt(\mathcal N ) $.

The purity of the exterior of a necklace is a consequence of the following
main result of the paper.\medskip

{\bf Theorem 3.} {\em  Let $\mathcal N $ be a Grassmann  necklace in
${[n]\choose r}$, $X\in \mathcal{O}ut(\mathcal N )$, and  $Y\in
\mathcal{I}nt(\mathcal N)$. Then $X$ and $Y$ are separated,
$X\|Y$. }\medskip

We prove this theorem in the next section. Now we establish its important
corollary.\medskip

{\bf Proposition 3.} {\em  Let $\mathcal N $ be a neklace. Then the exterior
$\mathcal{O}ut(\mathcal N )$ of $\mathcal N$ is a pure system; its rank is
equal to the number of alignments of the corresponding permutation
$\pi(\mathcal N)$.}\medskip

{\em Proof of Proposition 3}. Let $\mathcal C $ be a maximal separated
system in $\mathcal{O}ut(\mathcal N )$ and let $\mathcal D$ be a
maximal separated system in $\mathcal{I}nt(\mathcal N )$.
Obviously, $\mathcal N \subset \mathcal{D}$.

We claim that the union $\mathcal C \cup \mathcal D $ is a  maximal
separated system in the Grassmanian  ${[n]\choose r}$. Indeed, due to
Theorem 3, the union is separated. To see the maximality, suppose that the
union can be extended by an additional  set $Z$ of cardinality $r$. Since
$Z$ is separated from $\mathcal N $, it belongs to $\mathcal S(\mathcal N)
$. Hence $Z$ belongs  either to $\mathcal{I}nt(\mathcal N)$ or to
$\mathcal{O}ut(\mathcal N)$, which contradicts the maximality of $\mathcal
D $ or $\mathcal C $.

By Theorem 1, the size of $\mathcal C \cup \mathcal D $ does not depend of
a choice of $\mathcal C $ and $\mathcal D $ (implying the same property
for each of $\mathcal C $ and $\mathcal D $). This proves the purity of
$\mathcal{I}nt(\mathcal N)$  and $\mathcal{O}ut(\mathcal N)$. The
assertion on the rank of $\mathcal{O}ut(\mathcal N)$ follows from the fact
that the rank of $\mathcal{I}nt(\mathcal{N})$ is equal to $k(n-k)+1$ minus
the number of alignments for $\pi$. \hfill $\Box$\medskip

{\bf Remark 3}. It may seem that the above reasonings lead to a new proof of
the purity of $\mathcal{I}nt$. However, they rely on Theorem~3, and the proof
of the latter given in Section~6 uses arguments from~\cite{OPS}.\medskip

Proposition 3 can be generalized for the case of two (or more) necklaces. To
formulate such generalizations, we use a shorter notation. Namely, considering
two necklaces $\mathcal{N}_1,\mathcal{N}_2$, we will write $\mathcal{I}_k$ for
$\mathcal{I}nt(\mathcal{N}_k)$, and write $\mathcal{O}_k$ for
$\mathcal{O}ut(\mathcal{N}_k)$, $k=1,2$.\medskip

{\bf Proposition 4.} {\em Suppose that necklaces $\mathcal N _1$ and
$\mathcal N _2$ are separated from each other. Then the following four
systems are pure: $\mathcal{I}_1\cap \mathcal{I}_2$, $\mathcal{I}_1\cap
\mathcal{O}_2$, $\mathcal{O}_1\cap \mathcal{I}_2$, and $\mathcal{O}_1\cap
\mathcal{O}_2$. The sum of their ranks is equal to $r(n-r)+1$.}\medskip

\unitlength=1mm \special{em:linewidth 0.4pt} \linethickness{0.4pt}
\begin{picture}(110.00,37.00)(0,4)
\put(55.00,17.50){\oval(50.00,15.00)}
\put(87.50,17.50){\oval(45.00,19.00)}
\put(52.00,6.00){\makebox(0,0)[cc]{$\mathcal{N}_1$}}
\put(88.00,4.00){\makebox(0,0)[cc]{$\mathcal{N}_2$}}
\put(72.00,17.00){\makebox(0,0)[cc]{$\mathcal{I}_1\cap \mathcal{I}_2$}}
\put(48.00,17.00){\makebox(0,0)[cc]{$\mathcal{I}_1\cap \mathcal{O}_2$}}
\put(93.00,17.00){\makebox(0,0)[cc]{$\mathcal{O}_1\cap \mathcal{I}_2$}}
\put(72.00,33.00){\makebox(0,0)[cc]{$\mathcal{O}_1\cap \mathcal{O}_2$}}
\end{picture}
  \begin{center}
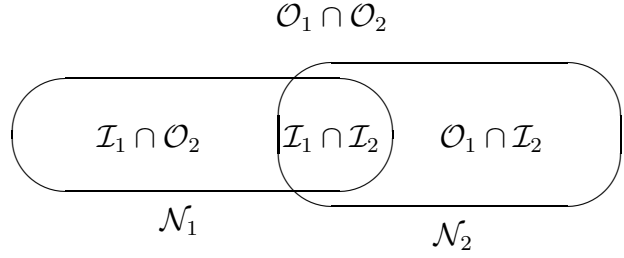
Figure 1. {\small Necklaces from Proposition 4.}
\end{center}

{\em Proof}. Let $\mathcal{A}$ be a maximal separated system in
$\mathcal{I}_1\cap \mathcal{I}_2$. Let $X\in \mathcal{N}_1\cap
\mathcal{I}_2$. Then $X\in\mathcal{I}_1\cap \mathcal{I}_2$,  ~$X$ is
separated from $\mathcal{I}_1$ and, moreover, $X$ is separated from
$\mathcal{A}$. By the maximality of $\mathcal{A}$, $X$ belongs to
$\mathcal{A}$. Therefore,\smallskip

a) $\mathcal{N}_1\cap \mathcal{I}_2$ (as well as $\mathcal{I}_1\cap
\mathcal{N}_2$) is contained in $\mathcal{A}$. \smallskip

Let $\mathcal{B}_1$ be a maximal separated system in $\mathcal{I}_1\cap
\mathcal{O}_2$. By similar reasonings, \smallskip

b1) $\mathcal{N}_1\cap \mathcal{O}_2$ is contained in $\mathcal{B}_1$.
\smallskip

Similarly, if $B_2$ is a maximal separated system in $O_1\cap I_2$, then
\smallskip

b2) $\mathcal{N}_2\cap \mathcal{O}_1$ is contained in $\mathcal{B}_2$.
\smallskip

Finally, let $\mathcal{C}$ be a maximal separated system in $\mathcal{O}_1\cap
\mathcal{O}_2$. We assert that the union $\mathcal{A}\cup \mathcal{B}_1\cup
\mathcal{B}_2\cup \mathcal{C}$ is a maximal separated system in  ${[n]\choose
r}$. Indeed:

First, by Theorem 3, this union is a  separated system.

Second, since $\mathcal{N}_1$ is separated from $\mathcal{N}_2$, we have
$\mathcal{N}_1=(\mathcal{N}_1\cap \mathcal{I}_2)\cup (\mathcal{N}_1\cap
\mathcal{O}_2)$. Hence, due to a) and b1), $\mathcal{N}_1$ is contained in
$\mathcal{A}\cup \mathcal{B}_1$. Similarly, $\mathcal{N}_2$ is contained in
$\mathcal{A}\cup \mathcal{B}_2$. Therefore, $\mathcal{N}_1$ and $\mathcal{N}_2$
are contained in $\mathcal{A}\cup \mathcal{B}_1\cup \mathcal{B}_2\cup
\mathcal{C}$.

Third,  let a set $X$ be separated from $\mathcal{A}\cup \mathcal{B}_1\cup
\mathcal{B}_2\cup \mathcal{C}$. Since $\mathcal{N}_1$ and $\mathcal{N}_2$ are
contained in the union $\mathcal{A}\cup \mathcal{B}_1\cup \mathcal{B}_2\cup
\mathcal{C}$, the set $X$ is separated from $\mathcal{N}_1$ and from
$\mathcal{N}_2$. Hence, $X$ belongs to one of the systems $\mathcal{I}_1\cap
\mathcal{I}_2$, ~$\mathcal{I}_1\cap \mathcal{O}_2$, ~$\mathcal{O}_1\cap
\mathcal{I}_2$, ~$\mathcal{O}_1\cap \mathcal{O}_2$. If $X$ belongs to
$\mathcal{I}_1\cap \mathcal{I}_2$, then it is separated from $\mathcal{A}$. By
the maximality of $\mathcal{A}$ in $\mathcal{I}_1\cap \mathcal{I}_2$, ~$X$
belongs to $\mathcal{A}$. In a similar way, we obtain $X\in\mathcal{A}$ in the
other cases. Thus, the maximality of the union is proven.

Now by Theorem 1, the size of the union $\mathcal{A}\cup \mathcal{B}_1\cup
\mathcal{B}_2\cup \mathcal{C}$ does not depend on the choice of
$\mathcal{A}$ in $\mathcal{I}_1\cap \mathcal{I}_2$. This proves the purity
of $\mathcal{I}_1\cap \mathcal{I}_2$. Similarly, we obtain the purity for
the other cases. \hfill$\Box$\medskip

There are two interesting special cases of necklaces $\mathcal N
_1,\mathcal N _2$ in Proposition 4. The first case is when one necklace is
`less' than the other.\medskip

{\bf Definition.} We say that $\mathcal{N}_1$ is {\em less} than
$\mathcal{N}_2$ if $\mathcal{I}nt(\mathcal{N}_1)\subset
\mathcal{I}nt(\mathcal{N}_2)$.\medskip

In this case, obviously, $\mathcal{N}_1\| \mathcal{N}_2$ and
$\mathcal{O}_2\subset \mathcal{O}_1$. We have the following
criterion:\medskip

{\bf Lemma 3.} A necklace {\em $\mathcal{N}_1$ is less than
a necklace $\mathcal{N}_2$ if and only if $\mathcal{N}_1\subset
\mathcal{I}nt(\mathcal{N}_2)$.}\medskip

{\em Proof}. The part `only if' is trivial because $\mathcal{N}_1\subset
\mathcal{I}nt(\mathcal{N}_1)$. Let us prove the converse: if
$\mathcal{N}_1\subset \mathcal{I}nt(\mathcal{N}_2)$, then
$\mathcal{I}nt(\mathcal{N}_1)\subset \mathcal{I}nt(\mathcal{N}_2)$.

Let ${X}\in \mathcal{I}nt(\mathcal{N}_1)$. We have to show that
$(\mathcal{N}_2)_i \ll_i X$ for any $i$, where $(\mathcal{N}_2)_i$ denotes
i-th set of the necklace $\mathcal{N}_2$. Since $(\mathcal{N}_2)_i \|
\mathcal{N}_1$ holds for any $i$, the set $(\mathcal{N}_2)_i$ belongs
either to $\mathcal{I}nt(\mathcal{N}_1)$ or to
$\mathcal{O}ut(\mathcal{N}_1)$.

In the first case, we have $(\mathcal{N}_1)_i\ll_i (\mathcal{N}_2)_i$ and
$(\mathcal{N}_2)_i\ll_i (\mathcal{N}_1)_i$, implying $(\mathcal{N}_1)_i=
(\mathcal{N}_2)_i$. Hence $(\mathcal{N}_2)_i= (\mathcal{N}_1)_i\ll_i X$.

In the second case,  $(\mathcal{N}_2)_i$ belongs to
$\mathcal{O}ut(\mathcal{N}_1)$. Then, by Theorem 3, $(\mathcal{N}_2)_i$ is
separated from $X$. Moreover, it holds that $(\mathcal{N}_2)_i\ll_i
(\mathcal{N}_1)_i$ (because $(\mathcal{N}_1)_i$ belongs to
$\mathcal{I}nt(\mathcal{N}_2)$) and $(\mathcal{N}_1)_i\ll_i X$ (because
$X$ belongs to $Int(\mathcal{N}_1)$). Thus, due to Lemma 1, we obtain
$(\mathcal{N}_2)_i\ll_i X$. \hfill$\Box$\medskip

The first special case is exposed in the following \medskip

{\bf Corollary 1.} {\em Let $\mathcal N_1 $ and $\mathcal N_2 $ be two
necklaces. Suppose that $\mathcal N_1 $ is less than  $\mathcal N_2 $,
$\mathcal{N}_1\subset \mathcal{I}_2$. Then the system $\mathcal{I}_2\cap
\mathcal{O}_1$ (the `ring' between $\mathcal N_2 $ and $\mathcal N_1 $) is
pure and its rank is equal to
$r(n-r)+1-rk(\mathcal{I}_1)-rk(\mathcal{O}_2)$.}\medskip

\unitlength=1mm \special{em:linewidth 0.4pt} \linethickness{0.4pt}
\begin{picture}(96.00,31.00)(0,0)
\put(70.00,15.00){\oval(30.00,10.00)[]} \put(70.00,15.00){\oval(50.00,20.00)[]}
\put(86.00,10.00){\makebox(0,0)[cc]{$\mathcal{N}_1$}}
\put(96.00,6.00){\makebox(0,0)[cc]{$\mathcal{N}_2$}}
\put(70.00,15.00){\makebox(0,0)[cc]{$\mathcal{I}_1$}}
\put(70.00,22.00){\makebox(0,0)[cc]{$\mathcal{O}_1\cap \mathcal{I}_2$}}
\put(70.00,28.00){\makebox(0,0)[cc]{$\mathcal{O}_2$}}
\end{picture}

{\em Proof}.  By Lemma 3, $\mathcal{I}_1\subset \mathcal{I}_2$ and
$\mathcal{O}_1\supset \mathcal{O}_2$. Therefore, $I_1\cap I_2=I_1$,
$I_1\cap O_2=\emptyset$, and $O_1\cap O_2=O_2$. Since
$\mathcal{N}_1\subset \mathcal{I}_2$, we have $\mathcal{N}_1 \|
\mathcal{N}_2$. Now the result follows from Proposition 4.
\hfill$\Box$\medskip

The second special case strengthens the condition $\mathcal{N}_1\|
\mathcal{N}_2$.\medskip

{\bf Corollary 2.} {\em Let $\mathcal{N}_1$ and $\mathcal{N}_2$ be
two necklaces. Suppose that $\mathcal{I}_1\| \mathcal{N}_2$ and
$\mathcal{I}_2\| \mathcal{N}_1$. Then $\mathcal{I}_1\cup
\mathcal{I}_2$ is a  pure system. }\medskip

{\em Proof}. The condition $\mathcal{I}_1\| \mathcal{N}_2$ (or
$\mathcal{I}_2\| \mathcal{N}_1$) implies $\mathcal{N}_1\| \mathcal{N}_2$.
Thus, we can apply Proposition 4. Moreover, the relation  $\mathcal{I}_1\|
\mathcal{N}_2$ gives the partitions $\mathcal{I}_1=(\mathcal{I}_1\cap
\mathcal{I}_2)\sqcup (\mathcal{I}_1\cap \mathcal{O}_2)$ and
$\mathcal{I}_2=(\mathcal{I}_2\cap \mathcal{I}_1)\sqcup (\mathcal{I}_2\cap
\mathcal{O}_1)$. Therefore, we have the partition
 $$
\mathcal{I}_1\cup \mathcal{I}_2=(\mathcal{I}_1\cap \mathcal{I}_2)\sqcup
(\mathcal{I}_1\cap \mathcal{O}_2) \sqcup (\mathcal{I}_2\cap
\mathcal{O}_1).
 $$
Let $\mathcal{C}$ be a maximal separated system in $\mathcal{I}_1\cup
\mathcal{I}_2$. Consider the intersection of $\mathcal{C}$ with each of
$\mathcal{I}_1\cap \mathcal{I}_2$, ~$\mathcal{I}_1\cap \mathcal{O}_2$, and
$\mathcal{I}_2\cap \mathcal{O}_1$. We assert that $\mathcal{C}\cap
(\mathcal{I}_1\cap \mathcal{I}_2)$ is a maximal separated system in
$\mathcal{I}_1\cap \mathcal{I}_2$. Indeed, suppose that one can extend it
by adding a new set $X\in \mathcal{I}_1\cap \mathcal{I}_2$. Since $X$ is
separated from $\mathcal{I}_1\cap \mathcal{O}_2$ and from
$\mathcal{I}_2\cap \mathcal{O}_1$, ~$X$ is separated from $\mathcal{C}$. A
contradiction. Similarly, $\mathcal{C}\cap (\mathcal{I}_1\cap
\mathcal{O}_2)$ is maximal in $\mathcal{I}_1\cap \mathcal{O}_2$, and
$\mathcal{C}\cap (\mathcal{I}_2\cap \mathcal{O}_1)$ is maximal in
$\mathcal{I}_2\cap \mathcal{O}_1$. By Proposition 4,
$|\mathcal{C}|=rk(\mathcal{I}_1\cap \mathcal{I}_2)+rk(\mathcal{I}_1\cap
\mathcal{O}_2)+rk(\mathcal{I}_2\cap \mathcal{O}_1)$. \hfill$\Box$\medskip

Note that if, in addition to the hypotheses in Corollary 2, we require
that $\mathcal{I}_1$ and $\mathcal{I}_2$ are disjoint, then it follows
that $rk(\mathcal{I}_1\cup
\mathcal{I}_2)=rk(\mathcal{I}_1)+rk(\mathcal{I}_2)$.\medskip

\section{Proof of Theorem 3}

Theorem 3 can be reformulated in the following equivalent form.\medskip

{\bf Theorem $3'$.} {\em Let $\mathcal{N}$  be a necklace. Suppose
that $\mathcal{C}$ is  a maximal separated system in the
Grassmanian ${[n]\choose  r}$, containing $\mathcal{N}$, and
$\mathcal{C}'=\mathcal{C}\cap \mathcal{I}nt(N)$. Then
$\mathcal{C}'$ is a maximal separated system in
$\mathcal{I}nt(N)$.}\medskip

Indeed, let $X$ be a set in $\mathcal{I}nt(N)$ which is separated from
$\mathcal{C}'$. Then, due to Theorem~3, $X$ is separated from
$\mathcal{C}-\mathcal{C}'$. Therefore, $X$ is separated from $\mathcal{C}$. By
the maximality of $\mathcal{C}$, ~$X$ belongs to $\mathcal{C}$ and, hence,
belongs to $\mathcal{C}'$.

To prove the converse, we notice that Theorem $3'$ can be regarded as a
generalization of the following\medskip

{\bf Theorem 4.} \cite[Theorem 3]{Sc}, see also \cite[Proposition 3.2]{OPS}
{\em Let $A$ be a subset in $[n]\choose {r-2}$, and let $i,j,k,l$ be a
cyclically ordered quadruple of elements of $[n]-A$. Suppose that $\mathcal{C}$
is a maximal separated system in the Grassmanian $[n]\choose {r}$ containing
the sets $Aij$, $Ajk$, $Akl$, $Ali$. Then $\mathcal{C}$ contains either $Aik$
or $Ajl$.}\medskip

Here we can interpret the quadruple $Aij$, $Ajk$, $Akl$, $Ali$ as a `small
necklace' whose interior consists of the quadruple plus the sets $Aik$ and
$Ajl$. There are two maximal separated systems in the interior of this
necklace, one containing $Aik$ and the other containing $Ajl$. Moving from one
of such systems to the other is called a {\em mutation}.\medskip

Let us deduce Theorem 3 from Theorem~$3'$. Let $\mathcal{N},X,Y$ be as in the
hypotheses of Theorem~3. Consider a maximal separated system $\mathcal{C}$ in
the Grassmannian containing $X$ and $\mathcal{N}$. Due to Theorem $3'$, its
restriction $\mathcal{C}'=\mathcal{C}\cap\mathcal{I}nt(\mathcal N)$ is a
maximal separated system in $\mathcal{I}nt(\mathcal{N})$. Let $\mathcal C ''$
be a maximal separated system in $\mathcal{I}nt(\mathcal N )$ which contains
$Y$. Due to Postnikov's theorem (\cite[Theorem~13.4]{Po}, see also
\cite[theorem 4.7]{OPS}), the systems $\mathcal C '$ and $\mathcal C ''$ can be
connected by a sequence of mutations. Each mutation preserves the separation
from $X$ (Theorem 4). Therefore, $X$ is separated from $\mathcal C ''$, and we
get $X\| Y$. \hfill$\Box$\medskip

Thus, it remains to prove Theorem~$3'$. Using a decomposition of the necklace
along with the corresponding permutation and the interior of the necklace into
connected components \cite[Sec.~5]{OPS}, one may assume that the necklace
$\mathcal{N}$ is {\em connected}, that is the sets $N_i$, $i\in [n]$, are
distinct. The proof will use a technique of plabic tilings developed in
\cite[Sec.~9]{OPS}. Let us recall this notion and details.\medskip

{\em Plabic tilings.} Suppose that $\mathcal C $ is a separated system in the
Grassmanian ${[n]\choose r}$. Then it is possible to construct a planar
bicolored (plabic) polygonal complex $\Sigma (\mathcal C )$, with a chessboard
coloring of its two-dimensional cells. In the beginning, we take $n$ vectors
$\xi _1,\ldots,\xi _n$ in the plane $\mathbb R ^2$, being the clockwise ordered
roots of $1$ of degree $n$ (identifying the plane with the set $\mathbb C$ of
complex numbers).

Then one can assign to every set $X\subset [n]$ the vector (point) $\xi (X)=
\sum_{i\in X}\xi _i$.

The set (structure) $\Sigma (\mathcal C)$ consists of 0-dimensional cells
(points), 1-dimensional cells (edges) and 2-dimensional cells (polygons), which
form a polygonal complex (where the nonempty intersection of two cells is again
a cell and is the common face of these two cells).

Here the 0-dimensional cells (vertices) are the points of the form $\xi (X)$
for $X\in \mathcal C $. One can check (using the separability) that these
points are distinct.

Two-dimensional cells are colored black and white. More precisely, let  $K$ be
an (r-1)-element subset of $[n]$. The {\em white clique}  $\mathcal W
(K)=\mathcal W _\mathcal C (K)$ consists of those sets $X\in \mathcal C $ that
contain $K$, $K\subset X$. Thus, $\mathcal W (K)$ consists of sets $Ka_1, Ka_2,
\ldots, Ka_k$, where the elements $a_1,\ldots,a_k$ are taken in cyclic order. A
white clique is {\em nontrivial} if it has at least three elements. For a
nontrivial white clique $W(K)$, the convex hull of the points $\xi (X)$, $X\in
\mathcal W (K)$, is a white-colored cell of the complex $\Sigma (\mathcal C )$.

Similarly, for a set $L$ of the size $r+1$,  the \emph{black clique} $\mathcal
B (L)$ is constituted from those sets $X\in \mathcal C $ that are contained in
$L$. A nontrivial black clique $\mathcal B(L)$ generates the black-colored
two-dimensional cell to be the convex hull of points  $\xi (X)$, where $X$ runs
over the elements of $\mathcal B (L)$.

The set of one-dimensional cells (edges) consists of the edges of its
two-dimensional cells and the segments joining vertices $\xi (X)$ and $\xi (Y)$
such that $\mathcal W (X\cap Y)=\mathcal B (X\cup Y)=\{X,Y\}$.

Let us notice that only neighbors can be joined by an edge in the complex
$\Sigma (\mathcal C )$, where sets  $X$ and $Y$ (of the same size) are called
{\em neighbors} if the symmetric difference $(X-Y)\cup (Y-X)$ consists of
exactly two elements.

The picture below illustrates the plabic tiling for a certain set-system; here
the sets of the system are indicated at the vertices and the letters on tiles
indicate their colors. (A more sophisticated example of plabic tilings is given
in \cite[Fig. 9]{OPS}.)

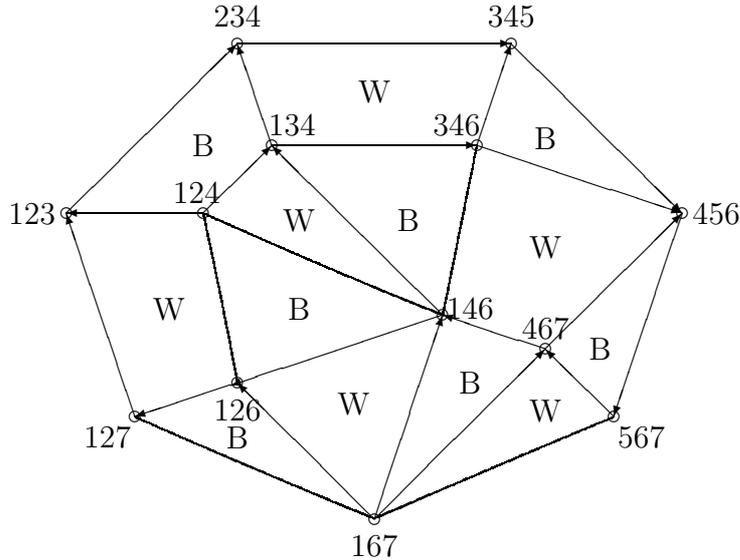
\unitlength=.9mm \special{em:linewidth 0.4pt}
\linethickness{0.4pt}
\begin{picture}(120.00,90.00)
\put(70.00,10.00){\circle{1.5}} \put(35.00,25.00){\circle{1.5}}
\put(105.00,25.00){\circle{1.5}} \put(95.00,35.00){\circle{1.5}}
\put(80.00,40.00){\circle{1.5}} \put(50.00,30.00){\circle{1.5}}
\put(45.00,55.00){\circle{1.5}} \put(25.00,55.00){\circle{1.5}}
\put(115.00,55.00){\circle{1.5}} \put(85.00,65.00){\circle{1.5}}
\put(55.00,65.00){\circle{1.5}} \put(50.00,80.00){\circle{1.5}}
\put(90.00,80.00){\circle{1.5}}
\bezier{200}(70.00,10.00)(86.00,17.00)(105.00,25.00)
\bezier{200}(70.00,10.00)(53.00,17.00)(35.00,25.00)
\put(35.00,25.00){\vector(-1,3){10.00}}
\put(45.00,55.00){\vector(-1,0){20.00}}
\put(25.00,55.00){\vector(1,1){25.00}}
\put(50.00,80.00){\vector(1,0){40.00}}
\put(90.00,80.00){\vector(1,-1){25.00}}
\put(115.00,55.00){\vector(-1,-3){10.00}}
\put(105.00,25.00){\vector(-1,1){10.00}}
\put(95.00,35.00){\vector(-3,1){15.00}}
\put(80.00,40.00){\vector(-3,-1){31.00}}
\put(49.00,29.67){\vector(-3,-1){14.00}}
\put(70.00,10.00){\vector(-1,1){20.00}}
\put(70.00,10.00){\vector(1,3){10.00}}
\put(80.00,40.00){\vector(-1,1){25.00}}
\put(55.00,65.00){\vector(1,0){30.00}}
\put(85.00,65.00){\vector(1,3){5.00}}
\put(45.00,55.00){\vector(1,1){10.00}}
\put(55.00,65.00){\vector(-1,3){5.00}}
\put(85.00,65.00){\vector(3,-1){30.00}}
\put(95.00,35.00){\vector(1,1){20.00}}
\bezier{150}(80.00,40.00)(82.00,50.00)(85.00,65.00)
\bezier{150}(45.00,55.00)(48.00,41.00)(50.00,30.00)
\bezier{200}(80.00,40.00)(65.00,46.00)(45.00,55.00)
\put(45.00,65.00){\makebox(0,0)[cc]{B}}
\put(70.00,73.00){\makebox(0,0)[cc]{W}}
\put(95.00,66.00){\makebox(0,0)[cc]{B}}
\put(95.00,50.00){\makebox(0,0)[cc]{W}}
\put(103.00,35.00){\makebox(0,0)[cc]{B}}
\put(70.00,10.00){\vector(1,1){25.00}}
\put(95.00,26.00){\makebox(0,0)[cc]{W}}
\put(84.00,30.00){\makebox(0,0)[cc]{B}}
\put(67.00,27.00){\makebox(0,0)[cc]{W}}
\put(50.00,22.00){\makebox(0,0)[cc]{B}}
\put(40.00,41.00){\makebox(0,0)[cc]{W}}
\put(59.00,41.00){\makebox(0,0)[cc]{B}}
\put(59.00,54.00){\makebox(0,0)[cc]{W}}
\put(75.00,54.00){\makebox(0,0)[cc]{B}}
\put(31.00,22.00){\makebox(0,0)[cc]{127}}
\put(20.00,55.00){\makebox(0,0)[cc]{123}}
\put(50.00,84.00){\makebox(0,0)[cc]{234}}
\put(90.00,84.00){\makebox(0,0)[cc]{345}}
\put(120.00,55.00){\makebox(0,0)[cc]{456}}
\put(109.00,22.00){\makebox(0,0)[cc]{567}}
\put(70.00,6.00){\makebox(0,0)[cc]{167}}
\put(50.00,26.00){\makebox(0,0)[cc]{126}}
\put(44.00,58.00){\makebox(0,0)[cc]{124}}
\put(58.00,68.00){\makebox(0,0)[cc]{134}}
\put(82.00,68.00){\makebox(0,0)[cc]{346}}
\put(95.00,38.00){\makebox(0,0)[cc]{467}}
\put(84.00,41.00){\makebox(0,0)[cc]{146}}
\end{picture}

Proposition  9.4 of  \cite{OPS} asserts that $\Sigma (\mathcal C )$ is a
complex.  In particular, the following is valid.\medskip

{\bf Fact.} {\em Let $X$ and $Y$ be neighbors of a separated system $\mathcal C
$. If the segment $[\xi (X), \xi      (Y)]$ and a cell  $C$ of $\Sigma(\mathcal
C)$ have more than one common point, then the points  $\xi (X)$ and $\xi (Y)$
are vertices of $C$.}\medskip

The tiling $\Sigma (\mathcal C )$ in the above picture fills is the regular
$n$-gon. This is not by chance, but is caused by the maximality of the system
$\mathcal C $.

Now let $\mathcal N =(N_1,\ldots,N_n)$ be a connected necklace. Let $\xi
(\mathcal N )$ be the closed polygonal curve (in the 1-dimensional subcomplex)
joining the points $\xi(N_1)$, $\xi(N_2),\ldots, ,\xi(N_n)$, $\xi(N_1)$ in this
order. An important fact (cf.~\cite[Proposition 8.8]{OPS}) is that $\xi
(\mathcal N )$ is a \emph{simple closed curve}. Therefore, it divides the plane
into the \emph{inside} and the \emph{outside} w.r.t. $\xi (\mathcal N )$, where
the former is homeomorphic to a disk and is denoted by $in(\mathcal N )$.

We reformulate Proposition 9.10 from \cite{OPS} as follows.\medskip

{\bf Proposition 5.} {\em Let $X\in {[n]\choose r}$ be separated from a
(connected) necklace $\mathcal N $. Then  $X$ belongs to
$\mathcal{I}nt(\mathcal N)$ if and only if $\xi (X)$ belongs to $in(\mathcal N
)$. } \hfill$\Box$\medskip

There is the following important characterization for the maximality of a
separated system established in \cite{OPS}.\medskip

{\bf Proposition 6.} {\em  Let $\mathcal N $ be a connected necklace, and let
$\mathcal C$ be a separated system in $\mathcal{I}nt(\mathcal N )$. The system
$\mathcal C $ is maximal in $\mathcal{I}nt(\mathcal N )$ if and only if the
complex $\Sigma (\mathcal C )$ fills in the polygon $in(\mathcal N )$.}\medskip

(One implication, namely, that the maximality of $\mathcal C $ implies
filling-in is stated in \cite[Proposition 11.2]{OPS}. For the converse
implication, let $\Sigma (\mathcal C)$ fill in $in(\mathcal N )$. Then the
graph $G$ dual to  $\Sigma (\mathcal C )$  is a reduced plabic graph (see the
proof of \cite[Proposition 11.2]{OPS}) and $\mathcal F (G)=\mathcal C $. Now
from \cite[Theorem 9.16]{OPS} it follows that $\mathcal F (G)$ is a maximal
separated system in $\mathcal{I}nt(\mathcal N )$.)\medskip

In particular, if $\mathcal C $ is a  separated system in
$\mathcal{I}nt(\mathcal N )$, then the complex $\Sigma (\mathcal C )$ is
located in the polygon $in(\mathcal N )$.

Now we are ready to prove the theorem.\medskip

{\em Proof of Theorem $3'$.} Let $\mathcal C $ be a maximal separated system in
${[n]\choose r}$. Then the complex $\Sigma (\mathcal C )$ fills in the regular
$n$-gon. Let $\mathcal N $ be a connected necklace, $\xi (\mathcal N )$ the
corresponding simple closed polygonal curve, and $in(\mathcal N )$ the inside
of the curve.

The intersection $\mathcal C '=\mathcal C \cap \mathcal{I}nt(\mathcal N )$ is a
separated system in $\mathcal{I}nt(\mathcal N )$.

\unitlength=1mm \special{em:linewidth 0.6pt} \linethickness{1pt}
\begin{picture}(120.00,90)(0,5)
\put(70.00,10.00){\circle{1.50}} \put(35.00,25.00){\circle{1.50}}
\put(105.00,25.00){\circle{1.50}} \put(95.00,35.00){\circle{1.50}}
\put(80.00,40.00){\circle{1.50}} \put(50.00,30.00){\circle{1.50}}
\put(45.00,55.00){\circle{1.50}} \put(25.00,55.00){\circle{1.50}}
\put(115.00,55.00){\circle{1.50}} \put(85.00,65.00){\circle{1.50}}
\put(55.00,65.00){\circle{1.50}} \put(50.00,80.00){\circle{1.50}}
\put(90.00,80.00){\circle{1.50}}
\bezier{200}(70.00,10.00)(86.00,17.00)(105.00,25.00)
\bezier{200}(70.00,10.00)(53.00,17.00)(35.00,25.00)
\put(35.00,25.00){\vector(-1,3){10.00}}
\put(45.00,55.00){\vector(-1,0){20.00}}
\put(25.00,55.00){\vector(1,1){25.00}}
\put(50.00,80.00){\vector(1,0){40.00}}
\put(90.00,80.00){\vector(1,-1){25.00}}
\put(115.00,55.00){\vector(-1,-3){10.00}}
\put(105.00,25.00){\vector(-1,1){10.00}}
\put(95.00,35.00){\vector(-3,1){15.00}}
\put(80.00,40.00){\vector(-3,-1){31.00}}
\put(49.00,29.67){\vector(-3,-1){14.00}}
\put(70.00,10.00){\vector(-1,1){20.00}}
\put(70.00,10.00){\vector(1,3){10.00}}
\put(80.00,40.00){\vector(-1,1){25.00}}
\put(55.00,65.00){\vector(1,0){30.00}}
\put(85.00,65.00){\vector(1,3){5.00}}
\put(45.00,55.00){\vector(1,1){10.00}}
\put(55.00,65.00){\vector(-1,3){5.00}}
\put(85.00,65.00){\vector(3,-1){30.00}}
\put(95.00,35.00){\vector(1,1){20.00}}
\bezier{150}(80.00,40.00)(82.00,50.00)(85.00,65.00)
\bezier{150}(45.00,55.00)(48.00,41.00)(50.00,30.00)
\bezier{200}(80.00,40.00)(65.00,46.00)(45.00,55.00)
\put(45.00,65.00){\makebox(0,0)[cc]{B}}
\put(70.00,73.00){\makebox(0,0)[cc]{W}}
\put(95.00,66.00){\makebox(0,0)[cc]{B}}
\put(95.00,50.00){\makebox(0,0)[cc]{W}}
\put(103.00,35.00){\makebox(0,0)[cc]{B}}
\put(70.00,10.00){\vector(1,1){25.00}}
\put(95.00,26.00){\makebox(0,0)[cc]{W}}
\put(84.00,30.00){\makebox(0,0)[cc]{B}}
\put(67.00,27.00){\makebox(0,0)[cc]{W}}
\put(50.00,22.00){\makebox(0,0)[cc]{B}}
\put(40.00,41.00){\makebox(0,0)[cc]{W}}
\put(59.00,41.00){\makebox(0,0)[cc]{B}}
\put(59.00,54.00){\makebox(0,0)[cc]{W}}
\put(75.00,54.00){\makebox(0,0)[cc]{B}}
\put(31.00,22.00){\makebox(0,0)[cc]{127}}
\put(20.00,55.00){\makebox(0,0)[cc]{123}}
\put(50.00,84.00){\makebox(0,0)[cc]{234}}
\put(90.00,84.00){\makebox(0,0)[cc]{345}}
\put(120.00,55.00){\makebox(0,0)[cc]{456}}
\put(109.00,22.00){\makebox(0,0)[cc]{567}}
\put(70.00,6.00){\makebox(0,0)[cc]{167}}
\put(50.00,26.00){\makebox(0,0)[cc]{126}}
\put(44.00,58.00){\makebox(0,0)[cc]{124}}
\put(58.00,68.00){\makebox(0,0)[cc]{134}}
\put(82.00,68.00){\makebox(0,0)[cc]{346}}
\put(95.00,38.00){\makebox(0,0)[cc]{467}}
\put(84.00,41.00){\makebox(0,0)[cc]{146}}
\bezier{52}(50.00,80.00)(68.00,73.00)(85.00,65.00)
\bezier{44}(85.00,65.00)(90.00,50.00)(95.00,36.00)
\bezier{24}(96.00,36.00)(101.00,31.00)(107.00,25.00)
\bezier{60}(107.00,25.00)(88.00,17.00)(70.00,9.00)
\bezier{56}(70.00,9.00)(51.00,17.00)(34.00,25.00)
\bezier{48}(34.00,25.00)(39.00,41.00)(44.00,55.00)
\bezier{44}(44.00,55.00)(47.00,68.00)(50.00,80.00)
\put(50.00,80.00){\circle*{2.00}}
\put(85.00,65.00){\circle*{2.00}}
\put(95.00,35.00){\circle*{2.00}}
\put(105.00,25.00){\circle*{2.00}}
\put(70.00,10.00){\circle*{2.00}}
\put(35.00,25.00){\circle*{2.00}}
\put(45.00,55.00){\circle*{2.00}}
\end{picture}

\begin{center}
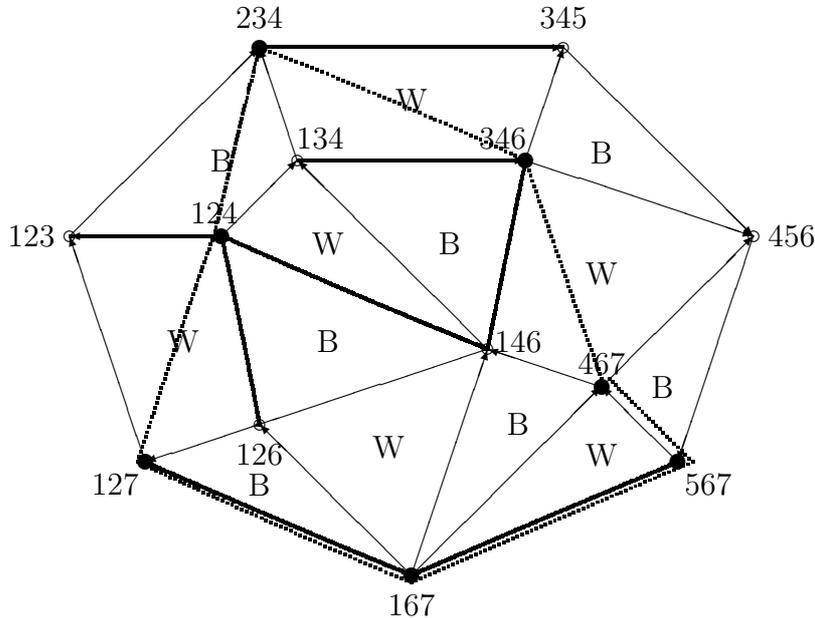
 Figure 2. \small{A necklace $\mathcal N $ is marked by bold dots. The curve
$\xi(\mathcal{N})$ is drawn by dotted lines.}
\end{center}

Hence, the complex  $\Sigma (\mathcal C ')$ is located inside the curve
$\xi(\mathcal N )$, ~$\Sigma(\mathcal{C})\subset in(\mathcal{N})$. We assert
that this complex fills in the polygon $in(\mathcal N )$. Indeed, let $P$ be a
point in  $in(\mathcal N )$. Since $\Sigma (\mathcal C )$ fills in the regular
$n$-gon, the point $P$ lies in some two-dimensional cell $C$ of $\Sigma
(\mathcal C )$. Let for definiteness $C$ be a white-colored cell corresponding
to a white clique $\mathcal W_\mathcal{C} (K)$. Consider the intersection of
$C$ with the polygon $in(\mathcal N )$. The edges of the polygonal curve $\xi
(\mathcal N )$ passing inside the cell  $C$ are some (non-intersecting)
diagonals of the convex polygon $C$ (see the Fact above and Fig. 3).

\unitlength=1mm \special{em:linewidth 0.6pt} \linethickness{1pt}
\begin{picture}(105.00,66.00)(-10,5)
\put(70.00,20.00){\vector(-2,1){10.00}}
\put(60.00,25.00){\vector(-1,2){5.00}}
\put(55.00,35.00){\vector(1,2){5.00}}
\put(60.00,45.00){\vector(2,1){10.00}}
\put(70.00,50.00){\vector(2,-1){10.00}}
\put(80.00,45.00){\vector(1,-2){5.00}}
\put(85.00,35.00){\vector(-1,-2){5.00}}
\put(80.00,25.00){\vector(-2,-1){10.00}}
\put(55.00,35.00){\vector(3,1){25.00}}
\put(80.00,25.00){\vector(-1,0){20.00}}
\bezier{34}(60.00,25.00)(20.00,30.00)(55.00,35.00)
\bezier{20}(81.00,44.00)(105.00,53.00)(81.00,59.00)
\bezier{28}(80.00,25.00)(105.00,20.00)(80.00,10.00)
\bezier{26}(80.00,10.00)(59.00,3.00)(40.00,10.00)
\bezier{26}(80.00,59.00)(60.00,66.00)(40.00,59.00)
\bezier{47}(40.00,59.00)(0.00,35.00)(40.00,10.00)
\put(70.00,34.00){\makebox(0,0)[cc]{$C$}}
\put(68.00,44.00){\makebox(0,0)[cc]{$P$}}

\put(72,44){\circle*{1.5}}

\put(38.00,44.00){\makebox(0,0)[cc]{$int(\mathcal{N})$}}
\end{picture}

\begin{center}
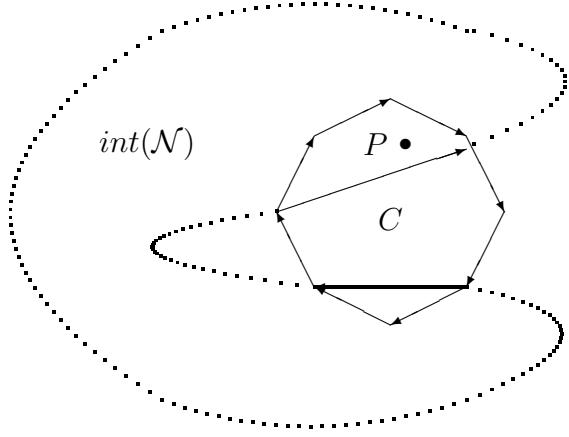
Figure 3. {\small Such a picture is impossible: the intersection of $C$ and
$in(\mathcal{N})$ is a convex polygon.}
\end{center}

\noindent Therefore, the intersection of $C$ with the polygon $in(\mathcal{N})$
is the union of convex polygons with vertices of the form $\xi(X)$, where $X\in
\mathcal{W}_\mathcal{C}(K)\cap \mathcal{C}'$. Hence, $P$ lies in the convex
hull of points $\xi (X)$, while $X$ runs over the set $\mathcal C '\cap
\mathcal W _{\mathcal C}(K)=\mathcal{W}_\mathcal{C'}(K)$. But then $P$ lies in
the white cell of the complex $\Sigma (\mathcal C ')$, corresponding to the
white clique       $\mathcal W _{\mathcal C '}(K)$.

Thus, the complex  $\Sigma (\mathcal C ')$ fills in the polygon $in(\mathcal N
)$ and, by Proposition 6, the separated system $\mathcal C '$ is maximal in
$\mathcal{I}nt(\mathcal N )$. Theorem $3'$ is proven. \hfill$\Box$\medskip

{\bf Remark 4}. Below we raise conjectures generalizing Theorem 3 (or $3'$),
Corollary 2, and some of results in~\cite{pur2}). Let
$\mathcal{K}=(K_1,\ldots,K_m)$ be a sequence of elements of the discrete
Grassmanian which satisfies the following three conditions:
\begin{enumerate}
\item $K_i$ and $K_{i+1}$ are neighbors for any $i=1,\ldots,m$ (where
$K_{m+1}=K_1$);
\item $\mathcal{K}$ is a separeted set-system;
\item the closed curve $\xi(\mathcal{K})$ is simple (without
self-intersections).
\end{enumerate}
We call such a $\mathcal{K}$ a {\em generalized necklace}. The inside
$in(\mathcal{K})$ of the curve $\xi(\mathcal{K})$ is defined as before. Define
the interior of the generalized necklace as follows:
 $$
\mathcal{I}nt(\mathcal{K}) =\{X\in {[n]\choose r}, X\| \mathcal{K}
\text{ and } \xi(X)\in in(\mathcal{K}) \}.
 $$
The exterior $\mathcal{O}ut(\mathcal{K})$ is defined to be the complement to
$\mathcal{I}nt(\mathcal{K})$ in $\mathcal{S}(\mathcal K)$:
~$\mathcal{O}ut(\mathcal{K})=\mathcal{S}(\mathcal K)\setminus
\mathcal{I}nt(\mathcal{K})$.
\medskip

{\bf Conjectures}.
\begin{enumerate}
\item {\em Both $\mathcal{I}nt(\mathcal{K})$ and
$\mathcal{O}ut(\mathcal{K})$) are pure systems.}
\item {\em If $C$ is a maximal separated system in the Grassmanian,
then the intersections $C\cap \mathcal{I}nt(K)$ and $C\cap
\mathcal{O}ut(K)$ are maximal separated systems in $ \mathcal{I}nt(K)$ and
$\mathcal{O}ut(K)$, respectively.}
\item $\mathcal{I}nt(K)\| \mathcal{O}ut(K)$.
\end{enumerate}

      \end{document}